\documentclass[a4paper, 12pt]{article}

\usepackage[utf8]{inputenc}
\usepackage[T1]{fontenc}
\usepackage[a4paper, margin=2.4cm]{geometry}
\usepackage{needspace}

\usepackage{libertine} 
\usepackage{inconsolata} 

\usepackage{amsmath, amsthm, amssymb}
\usepackage{graphicx}
\usepackage{enumerate}
\usepackage{authblk}
\usepackage[textsize=small, textwidth=2cm, color=yellow]{todonotes}
\usepackage{thmtools}
\usepackage{thm-restate}
\usepackage[colorlinks=true, citecolor=red]{hyperref}
\usepackage{float}

\declaretheorem[name=Theorem, numberwithin=section]{theorem}
\declaretheorem[name=Lemma, sibling=theorem]{lemma}

\declaretheorem[name=Corollary, sibling=theorem]{corollary}

\def\cqedsymbol{\ifmmode$\lrcorner$\else{\unskip\nobreak\hfil
\penalty50\hskip1em\null\nobreak\hfil$\lrcorner$
\parfillskip=0pt\finalhyphendemerits=0\endgraf}\fi}

\let\le\leqslant
\let\ge\geqslant

\title{Pivot-minors and the Erd\H{o}s-Hajnal conjecture}

\author[1]{James Davies}

\affil[1]{University of Cambridge, United Kingdom.}
\date{}

\begin{document}
	\sloppy

\maketitle

\begin{abstract}
    We prove a conjecture of Kim and Oum that every proper pivot-minor-closed class of graphs has the strong Erd\H{o}s-Hajnal property.
    More precisely, for every graph $H$, there exists $\epsilon > 0$ such that every $n$-vertex graph with no pivot-minor isomorphic to $H$ contains two sets $A, B$ of vertices such that $|A|, |B| \ge \epsilon n$ and $A$ is complete or anticomplete to $B$.
\end{abstract}

\section{Introduction}\label{EH:sec:intro}

For a graph $G$, we let $\alpha(G), \omega(G)$ denote the independence and clique number of $G$ respectively.
A class of graph $\mathcal{G}$ has the \emph{Erd{\H{o}}s-Hajnal property} if there exists a constant $c>0$ such that $\max\{\alpha(G),\omega(G)\}\ge |V(G)|^c$ for every $G\in \mathcal{G}$. Erd{\H{o}}s and Hajnal~\cite{erdos1989ramsey} conjectured that for every graph $H$, the graphs with no induced $H$ subgraph have the Erd{\H{o}}s-Hajnal property.

Strengthening the Erd{\H{o}}s-Hajnal property is the strong Erd{\H{o}}s-Hajnal property.
Two vertex sets $A,B\subseteq V(G)$ in a graph $G$ are \emph{complete} if $A,B$ are disjoint, and every vertex of $A$ is adjacent to every vertex of $B$.
Similarly, $A$ and $B$ are \emph{anticomplete} if they are disjoint and no vertex of $A$ is adjacent to a vertex of $B$.
A class of graphs $\mathcal{G}$ has the \emph{strong Erd{\H{o}}s-Hajnal property} if every graph $G\in \mathcal{G}$ with $|V(G)|>1$ contains two sets $A, B$ of vertices such that $|A|, |B| \ge \epsilon |V(G)|$ and $A$ is either complete or anticomplete to $B$.
The strong Erd{\H{o}}s-Hajnal property implies the Erd{\H{o}}s-Hajnal property, however many classes, such as triangle-free graphs do not satisfy the strong Erd{\H{o}}s-Hajnal property~\cite{fox2008erdHos}.

Chudnovsky and Oum~\cite{chudnovsky2018vertex} proved an analogue of the Erd{\H{o}}s-Hajnal conjecture that for every graph $H$, the class of $H$-vertex-minor-free graphs has the (strong) Erd{\H{o}}s-Hajnal property.
Kim and Oum~\cite{kim2021erdHos} conjectured a strengthening that for every graph $H$, the class of $H$-pivot-minor-free graphs has the (strong) Erd{\H{o}}s-Hajnal property. We prove this conjecture. Pivot-minors shall be formally defined in Section~\ref{EH:sec:pivot}.

\begin{theorem}\label{main}
    For every graph $H$, there exists $\epsilon > 0$ such that for all $n > 1$, every $n$-vertex graph not containing $H$ as a pivot-minor has two sets $A, B$ of vertices such
    that $|A|, |B| \ge \epsilon n$ and $A$ is either complete or anticomplete to $B$.
\end{theorem}

Since the strong Erd{\H{o}}s-Hajnal property implies the Erd{\H{o}}s-Hajnal property, we obtain the following corollary of Theorem~\ref{main}.

\begin{corollary}\label{EH:main}
    For every graph $H$, there exists $c > 0$ such that $\max\{\alpha(G),\omega(G)\} \ge |V(G)|^c$ for every graph $G$ not containing $H$ as a pivot-minor.
\end{corollary}

Theorem~\ref{main} was proved in the case that $H$ is a cycle by Kim and Oum~\cite{kim2021erdHos}.
Two other natural pivot-minor-closed classes are circle graphs and graphs of bounded rank-width. Both of these classes are also vertex-minor-closed, and so have the strong Erd{\H{o}}s-Hajnal property by the theorem of Chudnovsky and Oum~\cite{chudnovsky2018vertex}. The strong Erd{\H{o}}s-Hajnal property for circle graphs is also implied by an earlier theorem of Pach and Solymosi~\cite{pach2001crossing}. Almost tight bounds for the Erd{\H{o}}s-Hajnal property of circle graphs are now known. A theorem of the author~\cite{davies2022improved} implies that $\max\{\alpha(G),\omega(G)\} \ge \frac{1}{3}\sqrt{|V(G)|/\log_2 |V(G)|}$ for every circle graph $G$. Up to a constant factor, this is best possible due to a construction of Kostochka~\cite{kostochka2004coloring,kostochka1988upper}.

Recently the strong Erd{\H{o}}s-Hajnal property has been studied for other hereditary classes of graphs. For a graph $H$, we let $\overline{H}$ denote the \emph{complement} of $H$. Chudnovsky, Scott, Seymour, and Spirkl~\cite{chudnovsky2020pure} proved that for every pair of trees $T,F$, the class of graphs containing no induced $T$ or $\overline{F}$ subgraph has the strong Erd{\H{o}}s-Hajnal property.

A \emph{subdivision} of a graph $G$ is a graph obtainable from $G$ by replacing edges by possibly longer paths.
Chudnovsky, Scott, Seymour, and Spirkl~\cite{chudnovsky2021pure} also proved that for every graph $H$, the class of graphs containing no subdivision of $H$, or the complement of a subdivision of $H$ as an induced subgraph has the strong Erd{\H{o}}s-Hajnal property.
We shall consider a more general class of graphs.

A \emph{proper odd subdivision} of a graph $G$ is a graph obtainable from $G$ by replacing each edge by a path of odd length at least $3$.
A \emph{fuzzy odd path} is a graph obtained from an odd length path by possibly adding an additional edge $uv$ between two internal vertices of the path, such that $uv$ is contained in a triangle.
We say that a graph $G$ is a \emph{proper fuzzy odd subdivision} of a graph $H$ if $G$ can be obtained from $H$ by replacing each edge with a fuzzy odd path of length more than one.
To find an $r$-vertex graph $H$ as a pivot-minor, it is enough to instead find as an induced subgraph either a proper fuzzy odd subdivision of $K_r$ (see Lemma \ref{EH:lem:oddsub}) or the complement of a graph obtainable from $K_{r+1}$ by replacing each edge with a longer path of length $1 \pmod{3}$ (see Lemma \ref{lem:K^mod3v2}).
By finding induced proper fuzzy odd subdivisions, we actually also prove the following theorem, which generalizes the theorem of Chudnovsky, Scott, Seymour and, Spirkl~\cite{chudnovsky2021pure} on induced subdivisions and their complements.

\begin{theorem}\label{mainfuzzy}
    For every graph $H$, there exists $\epsilon > 0$ such that for all $n > 1$, every $n$-vertex graph containing no proper fuzzy odd subdivision of $H$ or the complement of a fuzzy odd subdivision of $H$ as an induced subgraph, contains two sets $A, B$ of vertices such
    that $|A|, |B| \ge \epsilon n$ and $A$ is either complete or anticomplete to $B$.
\end{theorem}

Our proof of Theorem~\ref{main} and Theorem~\ref{mainfuzzy} adapts and extends that of Chudnovsky, Scott, Seymour, and Spirkl~\cite{chudnovsky2021pure}.
It would be natural to conjecture a further strengthening of Theorem~\ref{mainfuzzy}. Namely, that Theorem~\ref{mainfuzzy} holds simply for odd subdivisions of $H$, rather than fuzzy odd subdivisions.
However, this turns out to be false. For sufficiently large $n$, Fox~\cite{fox2006bipartite} constructed $n$-vertex graphs containing no odd hole or antihole and with no complete or anticomplete pair $A, B\subseteq V(G)$ with
$|A|, |B| \ge \frac{15n}{\log_2 n}$.
However, Scott, Seymour, and Spirkl~\cite{scott2021pure} did prove that for every graph $H$ and $c>0$, there exists $\epsilon >0$ such that for all $n > 1$, every $n$-vertex graph containing no odd subdivision of $H$ or the complement of an odd subdivision of $H$ as an induced subgraph, has two sets $A, B$ of vertices such
that $|A| \ge \epsilon n$, $|B|\ge \epsilon n^{1-c}$, and $A$ is either complete or anticomplete to $B$.
If instead of forbidding odd subdivision, one instead forbids induced subdivisions where each edge is replaced with a longer path with length $2\ell \pmod{k}$, then these graphs do have the strong Erd{\H{o}}s-Hajnal property (see Theorem \ref{thm:EHmod}).

In Section~\ref{EH:sec:pivot} we formally introduce pivot-minor and prove all the necessary lemmas concerning them.
The rest of the paper will then focus on induced subgraphs. In Section~\ref{EH:sec:coherence} there is a further overview of the remainder of the proofs of Theorem~\ref{main} and Theorem~\ref{mainfuzzy}.

\section{Pivot-minors and subdivisions}\label{EH:sec:pivot}

For a vertex $v$ in a graph $G$, we let $N(v)$ denote its \emph{neighbourhood}, i.e., the set of vertices adjacent to $v$.
We let $N[v]=N(v)\cup \{v\}$. When the graph $G$ is not clear from context, we use $N_G(v), N_G[v]$.

For an edge $uv$ of a graph $G$, let $V_1 = N(u) \backslash N[v]$, $V_2 = N(v) \backslash N[u]$, and $V_3 = N(u) \cap N(v)$.
\emph{Pivoting} the edge $uv$ of $G$ is the act of first complementing the edges between each of the three pairs of vertex sets $(V_1, V_2)$, $(V_2, V_3)$, and $(V_1, V_3)$, and then swapping the vertex labels of $u$ and $v$.
We denote this graph by $G \wedge uv$.
A graph $H$ is a \emph{pivot-minor} of a graph $G$ if a graph isomorphic to $H$ can be obtained from $G$ by a sequence
of vertex deletions and pivots.

Given a graph $H$ and an integer $t\ge 1$, we denote by $H^t$ the subdivision of $H$ obtained by replacing each edge with a path of length $t+1$.
For a positive integer $r$, we let $[r]=\{1,\ldots , r\}$.
For a graph $G$ and $X\subseteq V(G)$, we let $G[X]$ denote the induced subgraph of $G$ on vertex set $X$.

We begin with a simple universal graph for pivot-minors. 

\begin{lemma}\label{EH:lem:K_r2}
The graph $K_r^2$ contains every graph on at most $r$ vertices as a pivot-minor.
\end{lemma}

\begin{proof}
Clearly it is enough to prove that $K_r^2$ contains every graph on exactly $r$ vertices as a pivot-minor.
So, let $H$ be a graph on vertex set $[r]$.
Let $V(K_r^2)= \{x_1,\ldots , x_r\} \cup \bigcup_{1\le i < j \le r} \{y_{i,j},z_{i,j}\}$, where for every pair $1\le i < j \le r$, $x_iy_{i,j}z_{i,j}x_j$ is a path.
Now, we find a pivot-minor isomorphic to $H$ by pivoting the edge $y_{i,j}z_{i,j}$ for every $ij\in E(H)$ with $i<j$, and then deleting the vertex set $\bigcup_{1\le i < j \le r} \{y_{i,j},z_{i,j}\}$.
\end{proof}

This can be extended to proper fuzzy odd subdivisions.

\begin{lemma}\label{EH:lem:oddsub}
Every proper fuzzy odd subdivision of $K_r$ contains every graph on at most $r$ vertices as a pivot-minor.
\end{lemma}

\begin{proof}
    Let $H$ be a proper fuzzy odd subdivision of $K_r$. If $H=K_r^2$, then the lemma follows by Lemma~\ref{EH:lem:K_r2}. So we argue inductively on $|V(H)|$. If $|V(H)|>|V(K_r^2)|$, then $H$ contains vertices $w,x,y,z$ such that either $wxyz$ is an induced path and $w,x,y,z$ all have degree 2 in $H$, or $N(x)=\{w,y,z\}$, $N(y)=\{x,z\}$, $w$ has degree 2 in $H$, and $z$ has degree 3 in $H$. In either case, $H'=(H \wedge xy) - \{x,y\}$ is a proper odd subdivision of $K_r^2$ with $|V(K_r^2)|\le |V(H')|=|V(H)|-2$. So, $H'$ (and thus $H$) contains every graph on at most $r$ vertices as a pivot-minor.
\end{proof}

For two graphs $G,H$, we let $G+H$ denote the graph obtained from the disjoint union of $G$ and $H$ by adding an edge between every vertex of $G$ and every vertex of $H$.

\begin{lemma}\label{lem:K^3}
    The graph $\overline{K_{r}^3}+K_1$ contains every graph on at most $r$ vertices as a pivot-minor.
\end{lemma}

\begin{proof}
    Clearly it is enough to prove that $\overline{K_{r}^3}+K_1$ contains every graph on exactly $r$ vertices as a pivot-minor.
    So, let $H$ be a graph on vertex set $[r]$.
    Let $V(\overline{K_{r}^3}+K_1)= \{x_1,\ldots , x_r\} \cup \bigcup_{1\le i < j \le r} \{x_{i,j}^1,x_{i,j}^2, x_{i,j}^3\} \cup \{u\}$, where for every pair $1\le i < j \le r$, $x_ix_{i,j}^1x_{i,j}^2x_{i,j}^3x_j$ is a path in $K_r^3$, and $u$ is the vertex of the $K_1$.
    Now, we find a pivot-minor isomorphic to $H$ by pivoting the edge $ux_{i,j}^2$ and then the edge $x_{i,j}^1x_{i,j}^3$ for every $ij\not\in E(H)$ with $i<j$, and then deleting the vertex set $\bigcup_{1\le i < j \le r} \{x_{i,j}^1,x_{i,j}^2, x_{i,j}^3\} \cup \{u\}$.
\end{proof}

\begin{lemma}\label{lem:K^mod3}
    Let $H$ be a graph obtained from $K_r$ by replacing each edge with a longer path of length $1 \pmod{3}$.
    Then, the graph $\overline{H} + K_1$ contains every graph on at most $r$ vertices as a pivot-minor.
\end{lemma}

\begin{proof}
    If $H= K_r^3$, then the result follows from Lemma \ref{lem:K^3}. So, we shall proceed inductively on $|V(H)|$.
    We may assume that $|V(H)|> |V(K_r^3)|$.
    Therefore, there exists an induced path $x_1x_2x_3x_4x_5$ of $H$ such that $x_2,x_3,x_4$ each have degree 2 in $H$ and $H'=(H \backslash \{x_2,x_3,x_4\}) \cup \{x_1x_5\}$ is a graph obtained from $K_r$ by replacing each edge with a longer path of length $1 \pmod{3}$.
    Since $|V(H')|<|V(H)|$, by inductive hypothesis, it follows that $\overline{H'} + K_1$ contains every graph on at most $r$ vertices as a pivot-minor.
    Therefore, it is enough to show that  $\overline{H} + K_1$ contains  $\overline{H'} + K_1$ as a pivot-minor.
    Let $u$ be the vertex of the $K_1$ in  $\overline{H} + K_1$.
    Then observe that  $((\overline{H} + K_1) \wedge ux_3 \wedge x_2x_4) \backslash \{x_2,x_3,x_4\} = \overline{H'} + K_1$, as desired.
\end{proof}

Immediately from this lemma, we obtain the following slightly more convenient lemma.

\begin{lemma}\label{lem:K^mod3v2}
    Let $H$ be a graph obtained from $K_{r+1}$ by replacing each edge with a longer path of length $1 \pmod{3}$.
    Then, the graph $\overline{H}$ contains every graph on at most $r$ vertices as a pivot-minor.
\end{lemma}

Let us remark that Lemma~\ref{lem:K^mod3v2} actually highlights one difference in the structure of pivot-minor-closed classes of graphs compared to vertex-minor-closed classes of graphs.
The complement $\overline{G}$ of a graph $G$ is a special case of a rank-$1$ perturbation of $G$.
Vertex-minors are well behaved under bounded rank perturbations. To be more precise, if $\mathcal{F}$ is a proper vertex-minor-closed class of graphs, then not every graph can be obtained as a vertex-minor of a bounded rank perturbation of a graph in $\mathcal{F}$ (see~\cite[Lemma 1.6.7]{rosethesis}).
However, this turns out to not be the case for pivot-minors.

Bipartite graphs are a proper pivot-minor-closed class of graphs and for each $r$, $K_{r+1}^3$ is bipartite, in particular $\overline{K_{r+1}^3}$ is a rank-1 perturbation of a bipartite graph.
Therefore, by Lemma~\ref{lem:K^mod3v2}, there exist proper pivot-minor-closed classes of graphs whose rank-1 perturbations contain all graphs as pivot-minors.
Another natural proper pivot-minor-closed class of graphs are the pivot-minors of line graphs~\cite{oum2009}.
As pointed out by Sang-il Oum, this class also contains $K_{r+1}^3$, thus this class also has the property that its rank-1 perturbations contain all graphs as pivot-minors.

We require another lemma on induced subgraphs that shall allow to more easily apply Lemma~\ref{lem:K^mod3v2}.

Let $G$ be a graph and let $F$ be a subgraph of $G$. Let $H$ be a graph obtained from $G$ by subdividing at least once every edge of $G$ not in $E(F)$, and not subdividing the edges in $E(F)$. We call such a graph $H$ (and graphs isomorphic to it) an \emph{$F$-filleting} of $G$.
Chudnovsky, Scott, Seymour, and Spirkl~\cite{chudnovsky2021pure} observed that for every pair of positive integers $\ell,k$, there exists a graph $G$ and a subgraph $F$ that is a forest whose components are all paths, such that every $F$-filleting of $G$ contains a hole of length $2\ell \pmod{k}$.
We shall require a slight extension for induced subdivisions rather than holes.

\begin{lemma}\label{lem:fillet}
    Let $\ell,k,r$ be a triple of positive integers. Then there exists a graph $G$ with a subgraph $F$ that is a forest whose components are all paths, such that every $F$-filleting of $G$ contains as an induced subgraph a graph $J$ that is obtainable from $K_r$ by replacing each edge with a longer path of length $2\ell \pmod{k}$.
\end{lemma}

\begin{proof}
    Let $R=R_{k^2}(r+3kr^2+3k-1)$ be the multicolour Ramsey number such that every $k^2$-edge-colouring of $K_{R}$ contains a monochromatic $K_{r+3kr^2}$.
    Now let $G$ be obtained from two disjoint copies of $K_{R}$ on vertex sets $\{x_1,\ldots , x_{R} \}$ and $\{y_1,\ldots , y_R \}$ by for each $1\le i \le R$ adding a path $F_i$ between $x_i$ and $y_i$ of length $\ell$. Let $F$ be the union of these paths.

    Let $H$ be some $F$-filleting of $G$.
    For each pair $1\le i < j \le R$, let $X_{i,j}$ be the path in $H$ corresponding to the edge $x_ix_j$ of $G$, and similarly, let $Y_{i,j}$ be the path in $H$ corresponding to the edge $y_iy_j$ of $G$.
    By Ramsey's theorem, there exists some $T\subset R$ with $|T| = r+3kr^2+3k-1$ and some pair $0\le \ell_x,\ell_y <k$ such that for every pair $i,j\in T$ with $i<j$, we have that $X_{i,j}$ is a path of length $\ell_x  \pmod{k}$ and $Y_{i,j}$ is a path of length $\ell_y \pmod{k}$.
    Without loss of generality, we may assume that $T=[r+3kr^2+3k-1]$.

    Now for each pair $1\le i < j \le r$, let $P_{i,j}^1 = X_{i,r+3k(ir+j)+1} \cup \bigcup_{t=1}^{k-1} X_{r+3k(ir+j)+t, r+3k(ir+j)+t+1}$,
    $P_{i,j}^2 = \bigcup_{t=0}^{k-1} Y_{r+3k(ir+j)+k+t, r+3k(ir+j)+k+t+1}$, 
    $P_{i,j}^3 = \bigcup_{t=0}^{k-2} X_{r+3k(ir+j)+2k+t, r+3k(ir+j)+2k+t+1}
    \cup X_{r+3k(ir+j)+3k-1, j}$, and let $P_{i,j} = P_{i,j}^1 \cup F_{r+3k(ir+j)+k}  \cup P_{i,j}^2 \cup  F_{r+3k(ir+j)+2k} \cup P_{i,j}^3$.
    Notice that each path $P_{i,j}$ has length $2\ell \pmod{k}$.
    Then, taking $J$ to be the union of all such paths $P_{i,j}$ provides the desired induced subgraph of $G$.
\end{proof}

By Lemma \ref{lem:K^mod3v2} and Lemma \ref{lem:fillet}, we obtain the following.

\begin{lemma}\label{lem:filletpivot}
    For every graph $H$, there exists a graph $G$ with a subgraph $F$ that is a forest whose components are all paths, such that the complement of every $F$-filleting of $G$ contains $H$ as a pivot-minor.
\end{lemma}

\section{Coherence}\label{EH:sec:coherence}

In this section we shall introduce some required preliminaries and show that Theorem~\ref{main} and Theorem~\ref{mainfuzzy} are implied by Theorem~\ref{mainfuzzymass}, which in turn is implied by two lemmas, which each of the following two sections are then dedicated to proving.
As in~\cite{chudnovsky2021pure}, our results apply more generally to massed graphs, which we introduce next.

A function $\mu:2^{V(G)}\to \mathbb{R}$ is a \emph{mass} on a graph $G$ if
\begin{itemize}
    \item $\mu(\emptyset)=0$ and $\mu(V(G))=1$,
    \item $\mu(X) \le \mu(Y)$ for all $X\subseteq Y \subseteq V(G)$, and
    \item $\mu(X\cup Y) \le \mu(X) + \mu(Y)$ for all $X,Y\subseteq V(G)$.
\end{itemize}
We call such a pair $(G,\mu)$ a \emph{massed graph}. The typical example of a mass on a graph $G$ is $\mu(X)=|X|/|V(G)|$. Another example is $\mu(X)=\chi(G[X])/\chi(G)$.
For a vertex $v$ of a massed graph $(G,\mu)$, we also denote $\mu(\{v\})$ by simply $\mu(v)$.

For $\epsilon>0$, let us say that a massed graph $(G,\mu)$ is \emph{$\epsilon$-coherent} if
\begin{itemize}
    \item $\mu(v) < \epsilon$ for every $v\in V(G)$,
    \item $\mu(N(v)) < \epsilon$ for every $v\in V(G)$, and
    \item $\min\{\mu(A),\mu(B)\} < \epsilon$ for every anticomplete pair $A,B\subseteq V(G)$.
\end{itemize}

We will prove the following.

\begin{theorem}\label{mainfuzzymass}
    For every graph $H$, there exists $\epsilon > 0$ such that every $\epsilon$-coherent massed graph $(G,\mu)$ contains a proper fuzzy odd subdivision of $H$ as an induced subgraph.
\end{theorem}

By Lemma~\ref{EH:lem:oddsub}, we also obtain the following corollary of Theorem~\ref{mainfuzzymass}.

\begin{corollary}\label{EH:masscol}
    For every graph $H$, there exists $\epsilon > 0$ such that every $\epsilon$-coherent massed graph $(G,\mu)$ contains a pivot-minor isomorphic to $H$.
\end{corollary}

The following is R{\"o}dl's~\cite{rodl1986universality} theorem.

\begin{theorem}[R{\"o}dl~\cite{rodl1986universality}]\label{rodl}
For every graph $H$ and $\epsilon>0$, there exists $\delta>0$ such that for every graph $G$ containing no induced $H$ subgraph, there exists $X\subseteq V(G)$ with $|X|\ge \delta|V(G)|$ such that in one of $G[X]$, $\overline{G[X]}$, every vertex in $X$ has degree less than $\epsilon |X|$.
\end{theorem}

With R{\"o}dl's~\cite{rodl1986universality} theorem, we can prove Theorem~\ref{mainfuzzy} assuming Theorem~\ref{mainfuzzymass}.
This is just the standard argument of applying R{\"o}dl's~\cite{rodl1986universality} theorem, which is often used to prove the strong Erd\H{o}s-Hajnal property by splitting into a sparse case and a dense case.

\begin{proof}[Proof of Theorem~\ref{mainfuzzy} assuming Theorem~\ref{mainfuzzymass}.]
    For our graph $H$, choose $\epsilon' >0$ satisfying Theorem~\ref{mainfuzzymass}.
    Let $H'$ be some proper fuzzy odd subdivision of $H$.
    Now choose $\delta >0$ satisfying Theorem~\ref{rodl} for $H',\epsilon'$. Now, let $\epsilon= \epsilon' \delta$.

    Let $G$ be a graph containing no proper fuzzy odd subdivision of $H$ or the complement of a proper fuzzy odd subdivision of $H$ as an induced subgraph, and with $|V(G)|>1$.
    It now follows from Theorem~\ref{rodl} that there exists $X\subseteq V(G)$ with $|X|\ge \delta|V(G)|$ such that in one of $G[X]$, $\overline{G[X]}$, every vertex in $X$ has degree less than $\epsilon' |X|$.
    If every vertex in $X$ has degree less than $\epsilon' |X|$ in $G[X]$, then let $G^*=G[X]$, otherwise let $G^*= \overline{G[X]}$.
    For each $Y\subseteq X$, let $\mu(Y)=\frac{|Y|}{|X|}$. Then $(G^*,\mu)$ is a massed graph.
    Since $G$ contains no proper fuzzy odd subdivision of $H$ or the complement of a proper fuzzy odd subdivision of $H$ as an induced subgraph, $G^*$ contains no proper fuzzy odd subdivision of $H$ as an induced subgraph. So by Theorem~\ref{mainfuzzymass}, $G^*$ is not $\epsilon'$-coherent.

    By definition, we have that $\mu(N_{G^*}(v)) = |N_{G^*}(v)|/|X| <  \epsilon'|X|/|X|= \epsilon'$ for each $v\in X$.
    If $\frac{1}{|X|} \ge \epsilon'$, then $|V(G)| \le \frac{1}{\epsilon'\delta}=\frac{1}{\epsilon}$, and so there exist distinct vertices $u,v\in V(G)$ with $A=\{u\}, B=\{v\}$ either complete or anticomplete, and $|A|=|B|=1\ge \epsilon|V(G)|$.
    So, we may assume that $\frac{1}{|X|} < \epsilon'$.
    Therefore, $\mu(v) = 1/|X| < \epsilon'$ for all $v\in X$.
    Since $G^*$ is not $\epsilon'$-coherent, it follows that there exists an anticomplete pair $A,B\subseteq X$ with $\mu(A),\mu(B) \ge \epsilon'$.
    Therefore, $A,B$ are either complete or anticomplete in $G$, and $|A|,|B|\ge \epsilon'|X|\ge \epsilon'\delta |V(G)|\ge \epsilon|V(G)|$, as desired.
\end{proof}

As seen from the above proof, R{\"o}dl's~\cite{rodl1986universality} theorem splits proving the strong Erd\H{o}s-Hajnal property into handling a sparse and a dense case. Corollary \ref{EH:masscol} handles the sparse case for pivot-minors. For the dense case we require the following theorem of Chudnovsky, Scott, Seymour, and Spirkl~\cite{chudnovsky2021pure}, which is a technical strengthening of their result that for every graph $H$, the class of graphs containing no subdivision of $H$, or the complement of a subdivision of $H$ as an induced subgraph has the strong Erd{\H{o}}s-Hajnal property.

\begin{theorem}[Chudnovsky, Scott, Seymour, and Spirkl {\cite[2.3]{chudnovsky2021pure}}]\label{EH:thm:fillet}
    For every graph $H$ and every subgraph $F$ that is a forest whose components are all paths, there exists $\epsilon > 0$ such that every $\epsilon$-coherent massed graph $(G,\mu)$ contains an $F$-filleting of $H$ as an induced subgraph.
\end{theorem}

By Lemma \ref{lem:filletpivot}, Theorem \ref{EH:thm:fillet} then handles the dense case for pivot-minors.
Thus, to prove Theorem \ref{main}, it now just remains to prove Theorem~\ref{mainfuzzymass} (as this implies Corollary \ref{EH:masscol} which handles the sparse case).

Let us remark that by Lemma~\ref{lem:fillet}, one can also obtain the following two theorems from Theorem~\ref{EH:thm:fillet} and then from R{\"o}dl's~\cite{rodl1986universality} theorem (in the same way as in the proof of Theorem~\ref{mainfuzzy} assuming Theorem~\ref{mainfuzzymass}).

\begin{theorem}
    For every graph $H$, there exists $\epsilon > 0$ such that every $\epsilon$-coherent massed graph $(G,\mu)$ contains as an induced subgraph a graph $J$ that is obtainable from $H$ by replacing each edge with a longer path of length $2\ell \pmod{k}$.
\end{theorem}

\begin{theorem}\label{thm:EHmod}
    For every graph $H$, there exists $\epsilon > 0$ such that for all $n > 1$, every $n$-vertex graph containing no graph obtainable from $H$ by replacing each edge with a longer path of length $2\ell \pmod{k}$ and then possibly complementing as an induced subgraph, contains two sets $A, B$ of vertices such
    that $|A|, |B| \ge \epsilon n$ and $A$ is either complete or anticomplete to $B$.    
\end{theorem}

For a vertex $v$ of a graph $G$, we let $N^r(v)$ denote the set of all vertices at distance exactly $r$ from $v$ in $G$. Similarly, we let $N^r[v]$ denote the set of all vertex at distance at most $r$ from $v$. When the graph $G$ is not clear from context, we use $N^r_G(v), N^r_G[v]$. For a set $X\subseteq V(G)$ of vertices, we let $N^r[X]= \bigcup_{v\in X} N^r[v]$.
When $r=1$, we simply use the notation $N^1[X]=N[X]$.

A massed graph $(G, \mu)$ is \emph{$(\delta,r)$-focused} if for every $Z \subseteq V(G)$ with $\mu(Z) \ge \delta$, there is a vertex $v\in Z$ with $\mu(N^r_{G[Z]}[v])\ge \mu(Z)/2$.

For an integer $r\ge 1$ and $\epsilon > 0$, we say that a massed graph $(G,\mu)$ is \emph{$(\epsilon,r)$-coherent} if
\begin{itemize}
    \item $\mu(N^r[v]) < \epsilon$ for every vertex $v$, and
    \item $\min\{\mu(A),\mu(B)\} < \epsilon$ for every two anticomplete sets of vertices $A,B\subseteq V(G)$.
\end{itemize}
Note that every $(\epsilon,1)$-coherent graph is $\epsilon$-coherent.

The proof of Theorem~\ref{mainfuzzymass} is split into the following two lemmas.

\begin{lemma}\label{EH:lem:rcoherent}
    For every graph $H$, there exist $\epsilon>0$ and an integer $r\ge 1$, such that every $(\epsilon,r)$-coherent massed graph $(G,\mu)$ contains an induced proper fuzzy odd subdivision of $H$.
\end{lemma}

\begin{lemma}\label{EH:lem:focussedcoherent}
    For every graph $H$ and integer $r\ge 1$, there exist $\epsilon , \delta >0$ such that every $(\delta,r)$-focused $\epsilon$-coherent massed graph $(G,\mu)$ contains an induced proper fuzzy odd subdivision of $H$.
\end{lemma}

Let us prove Theorem~\ref{mainfuzzymass} assuming Lemma~\ref{EH:lem:rcoherent} and Lemma~\ref{EH:lem:focussedcoherent}.
This general strategy works for any hereditary class of graphs.

\begin{proof}[Proof of Theorem~\ref{mainfuzzymass} assuming Lemma~\ref{EH:lem:rcoherent} and Lemma~\ref{EH:lem:focussedcoherent}.]
    By Lemma~\ref{EH:lem:rcoherent} there exist $\epsilon_r > 0$ and an integer $r \ge 1$, such that every
$(\epsilon_r,r)$-coherent massed graph $(G,\mu)$ contains an induced proper fuzzy odd subdivision of $H$.
Choose $\epsilon_{r+1}, \delta>0$ satisfying the hypotheses of Lemma~\ref{EH:lem:focussedcoherent} for $r+1$.
Notice that $\epsilon_r$ and $\epsilon_{r+1}$ both could have been chosen to be arbitrarily small.
So, by possibly choosing one or both of $\epsilon_r$ and $\epsilon_{r+1}$ to be smaller, we may assume that $\epsilon_r\delta=\epsilon_{r+1}$ and $\epsilon_{r+1} \le \delta / 2$.

Let $(G,\mu)$ be an $\epsilon_{r+1}$-coherent massed graph.
By Lemma~\ref{EH:lem:focussedcoherent}, we may assume that $(G,\mu)$ is not $(\delta, r+1)$-focused, since otherwise $G$ would contain an induced proper fuzzy odd subdivision of $H$.
So, there exists $Z \subseteq V(G)$ with $\mu(Z) \ge \delta$, such that $\mu(N^{r+1}_{G[Z]} [v] ) < \mu(Z)/2$ for each $v\in Z$.
Let $\mu'(X)=\mu(X)/\mu(Z)$ for all $X\subseteq Z$. Then $(G[Z], \mu')$ is a massed graph.

Suppose for sake of contradiction that $(G[Z], \mu')$ is not $(\epsilon_r, r)$-coherent.
Note that $\mu(X) \ge \delta \mu'(X)$ for all $X\subseteq Z$. In particular, if $\mu'(X)\ge \epsilon_r$, then $\mu(X)\ge \epsilon_r\delta=\epsilon_{r+1}$.
So since $(G,\mu)$ is $\epsilon_{r+1}$-coherent, $G[Z]$ contains no anticomplete pair $A,B\subseteq Z$ with $\mu'(A),\mu'(B) \ge \epsilon_{r}$.
Hence, there exists $v\in Z$ such that $\mu'(N_{G[Z]}^{r}[v]) \ge \epsilon_{r}$.
Therefore $\mu(N_{G[Z]}^{r}[v]) \ge \epsilon_{r}\delta = \epsilon_{r+1}$.
Since $\mu(N_{G[Z]}^{r+1}[v]) < \mu(Z)/2$, we have that $\mu(Z\backslash N_{G[Z]}^{r+1}[v]) > \mu(Z)/2 \ge \delta/2 \ge \epsilon_{r+1}$.
This contradicts the fact that $(G,\mu)$ is $\epsilon_{r+1}$-coherent since $N_{G[Z]}^{r}[v]$ and $Z\backslash N_{G[Z]}^{r+1}[v]$ are anticomplete.
Hence, $(G[Z], \mu')$ is $(\epsilon_r, r)$-coherent, and the result follows from Lemma~\ref{EH:lem:rcoherent}.
\end{proof}

So now proving Theorem~\ref{mainfuzzymass} (and thus Corollary~\ref{EH:masscol}, Theorem~\ref{main}, Corollary~\ref{EH:main}, and Theorem~\ref{mainfuzzy}) just relies on proving Lemma~\ref{EH:lem:rcoherent} and Lemma~\ref{EH:lem:focussedcoherent}.
We will prove Lemma~\ref{EH:lem:rcoherent} in Section~\ref{EH:sec:small}, and Lemma~\ref{EH:lem:focussedcoherent} in Section~\ref{EH:sec:focus}.
For both of these proofs, we shall need a lemma of Chudnovsky, Scott, Seymour, and Spirkl~\cite{chudnovsky2021pure}. The rest of this section is dedicated to introducing this lemma. We require some further definitions first.

Let $\mathcal{Y}$ be a set of pairwise disjoint subsets of $V(G)$, where $G$ is a graph.
Let $N$ be a graph, and for each $v \in V(N)$ let $X_v \subseteq V(G)$. We say that the
family $(X_v : v \in V(N))$ is \emph{$\mathcal{Y}$-spread} if for each $v \in V(N)$ there exists $Y_v \in \mathcal{Y}$ such that the sets $(Y_v : v \in V(N))$ are all different, and $X_v \subseteq Y_v$ for each
$v\in V(N)$.

For disjoint $X,Y \subseteq V(G)$, we say that $X$ \emph{covers} $Y$ if every vertex in $Y$ has a neighbour in $X$. 
Let $(G,\mu)$ be a massed graph. We say that $X \subseteq V(G)$ is \emph{$\delta$-dominant} if $\mu(N[X]) \ge \delta$.

A \emph{leaf} of a tree $T$ is a vertex of degree at most 1. We let $L(T)$ denote the set of leaves of a tree $T$.
A \emph{caterpillar} is a tree $T$ such that $T -L(T)$ is a path.
A \emph{rooted caterpillar} is a tree $T$ with a distinguished vertex $h$, called its
\emph{head}, such that $T -L(T)$ is a path with one end being $h$.

Let $N$ be a rooted caterpillar with head $h$.
Let $(G,\mu)$ be a massed graph.
A \emph{$\delta$-realization} of $N$ in $G$ is an assignment of a subset $X_v \subseteq V(G)$ to each vertex $v\in V(N)$, satisfying the following conditions:
\begin{itemize}
    \item the sets $X_v$ $(v\in V(N))$ are pairwise disjoint,
    \item for every edge $uv$ of $N$, if $v$ lies on the path of $N$ between $u$ and $h$, then $X_u$ covers $X_v$,
    \item for all distinct $u,v\in V(N)$, if $u,v$ are nonadjacent in $N$, then $X_u,X_v$ are anticomplete, and
    \item $\mu(X_h)\ge \delta$, and for each $u\in V(N)\backslash \{h\}$, $X_u$ is connected and $\delta$-dominant.
\end{itemize}

For $X \subseteq V(G)$ and $v \in X$, let us say $v$ is an \emph{$r$-centre} of $X$ if every vertex in $X$ has $G[X]$-distance at most $r$ from $v$. Note in particular that if $X \subseteq V(G)$ has an $r$-centre, then $G[X]$ is connected.

The following lemma is due to Chudnovsky, Scott, Seymour, and Spirkl~\cite{chudnovsky2021pure}.

\begin{lemma}[Chudnovsky, Scott, Seymour, and Spirkl {\cite[5.2]{chudnovsky2021pure}}]\label{EH:lem:realization}
    Let $N$ be a rooted caterpillar, let $\delta, \epsilon > 0$, let $(G,\mu)$ be an $\epsilon$-coherent massed graph, and let $\mathcal{Y}$ be a set of disjoint subsets of $V(G)$ such that $|\mathcal{Y}| = 2^{|V(N)|}$ and $\mu(Y) \ge 2^{2^{|V(N)|}}(\delta + \epsilon)$ for each $Y \in \mathcal{Y}$. Then there is a $\mathcal{Y}$-spread $\delta$-realization of $N$ in $G$.
    
    If in addition $r \ge 0$ is an integer, $\epsilon \le \delta/2$, and $(G,\mu)$ is $(\delta,r)$-focused, then there is a $\mathcal{Y}$-spread $\delta$-realization $(X_v : v \in V(N))$ of $N$ in $G$ such that $X_v$ has an $r$-centre for each $v \in V(N)$, except the head.
\end{lemma}

The case that $(G,\mu)$ is $(\delta,r)$-focused is called the focused case of Lemma~\ref{EH:lem:realization}, and the other case is called the unfocused case. The unfocused case will be used in Section~\ref{EH:sec:small} to prove Lemma~\ref{EH:lem:rcoherent}, and the focused case will be used in Section~\ref{EH:sec:focus} to prove Lemma~\ref{EH:lem:focussedcoherent}.

\section{When small balls have small mass}\label{EH:sec:small}

This section is dedicated to proving Lemma~\ref{EH:lem:rcoherent}.
We require another lemma of Chudnovsky, Scott, Seymour, and Spirkl~\cite{chudnovsky2021pure}, but first we must introduce ladders.

A \emph{$k$-ladder} in a graph $G$
is a family of $3k$ subsets
\[
A_1, \ldots , A_k, B_1, \ldots , B_k, C_1, \ldots , C_k
\]
of $V(G)$, pairwise disjoint and such that
\begin{itemize}
    \item for $1\le i \le k$, $A_i$ is connected and covers $B_i$, and $B_i$ covers $C_i$,
    \item for $1\le i \le k$, $A_i,C_i$ are anticomplete, and
    \item for all distinct $i, j \in \{1,\ldots, k\}$, $A_i$ is anticomplete to $A_j \cup B_j \cup C_j$.
\end{itemize}
A $k$-ladder is \emph{half-cleaned}, if $B_i$ is anticomplete to $C_j$ for every $1\le i < j \le k$.
The \emph{union} of a $k$-ladder is the triple $(A,B,C) = A \cup B \cup C$, where $A=\bigcup_{i=1}^kA_i$, $B=\bigcup_{i=1}^k B_i$, and $C=\bigcup_{i=1}^k C_i$.

Chudnovsky, Scott, Seymour, and Spirkl~\cite{chudnovsky2021pure} proved the following.

\begin{lemma}[Chudnovsky, Scott, Seymour, and Spirkl  {\cite[6.2]{chudnovsky2021pure}}]\label{EH:lem:ladder}
Let $\epsilon,\kappa >0$, and let $k\ge 0$ be an integer such that $(k-1)k(\kappa + \epsilon)\le 1$ and
$(k-1)(\kappa + \epsilon)+4\epsilon\le 1$. Let $(G,\mu)$ be an $\epsilon$-coherent massed graph. Then there
is a half-cleaned $k$-ladder
\[
A_1, \ldots , A_k, B_1, \ldots , B_k, C_1, \ldots , C_k
\]
in $G$ such that $\mu(C_i)\ge \kappa$ for $1\le i \le k$.
\end{lemma}

The following lemma is implicit in the proof of {\cite[6.4]{chudnovsky2021pure}} due to Chudnovsky, Scott, Seymour, and Spirkl, and we closely follow their proof.
This lemma is slightly stronger than we need, but the stronger statement better facilitates the proof.

\begin{lemma}\label{EH:lem:caterpillarladder}
Let $\epsilon,\kappa >0$, and let $N$ be a caterpillar with $k$ leaves, with the property that ${(2^{|V(N)|}-1)2^{|V(N)|}(2^{2^{|V(N)|}}(\kappa + 2\epsilon) + \epsilon)\le 1}$ and
$(2^{|V(N)|}-1)(2^{2^{|V(N)|}}(\kappa + 2\epsilon) + \epsilon)+4\epsilon\le 1$.
Let $(G,\mu)$ be an $(\epsilon,5\cdot 2^{|V(N)|})$-coherent massed graph.
Then $G$ has an induced subgraph $N'$ isomorhpic to a subdivision of $N$ with leaves $L'=\{b_1,\ldots, b_k \}$ and a half-cleaned $k$-ladder
\[
A_1, \ldots , A_k, B_1, \ldots , B_k, C_1, \ldots , C_k,
\]
such that
\begin{itemize}
    \item $\mu(C_i)\ge \kappa$ for $1\le i \le k$,
    \item $N'$ and $(A,B,C)$ are disjoint,
    \item $V(N')\backslash L'$ and $(A,B,C)$ are anticomplete,
    \item for each $i\in \{1,\ldots , k\}$, $N(b_i)\cap (A,B,C)\subseteq A_{i}$, and 
    \item for each $i\in \{1,\ldots , k\}$, $G[A_{i}\cup \{b_i\}]$ is connected.
\end{itemize}
\end{lemma}

\begin{proof}
    Assign a head to $N$ by choosing an end of the path obtained from $N$ by deleting all leaves of $N$.
    This makes $N$ rooted.
    Let $k' = 2^{|V(N)|}$.
    By Lemma~\ref{EH:lem:ladder}, there is a half-cleaned $k'$-ladder
\[
A_1, \ldots , A_{k'}, B_1, \ldots , B_{k'}, C_1, \ldots , C_{k'}
\]
in $G$ such that $\mu(C_i) \ge 2^{k'}(\kappa + 2\epsilon)$ for $1 \le i \le k'$. Then the unfocused
case of Lemma~\ref{EH:lem:realization} (taking $\delta = \kappa + \epsilon)$ implies that there is a $\{C_1,\ldots,C_{k'}\}$-spread $(\kappa + \epsilon)$-realization $(X_v : v\in V(N))$ of $N$ in $G$.

Let $n_1, \ldots , n_q$ be the vertices of $N$ that are not leaves, where $n_1$ is the head of $N$ and $n_in_{i+1}$ is an edge of $N$ for $1 \le i < q$. Choose $x_1 \in X_{n_1}$.
For $2 \le i \le q$ in order, since $X_{n_i}$ covers $X_{n_{i-1}}$ by definition of a realization, we may choose a vertex $x_i \in X_{n_i}$ adjacent to $x_{i-1}$. Since $\{x_1, \ldots ,x_q \} \subseteq C_1 \cup \cdots \cup C_{k'}$, it follows that $\{x_1, \ldots ,x_q \}$ is anticomplete to $A_1 \cup \cdots \cup A_{k'}$. Also, since there are no edges between $X_u$ and $X_v$ for nonadjacent $u,v\in V(N)$, it follows that $x_1, \ldots ,x_q$ are the vertices in order of an induced path of $G$.
Since $(X_v : v\in V(N))$ is a realization of $N$ in $G$, we have that for each leaf $\ell$ of $N$ with neighbour $n_j$ say, $x_j$ has a neighbour in $X_\ell$,
and $X_\ell$ is anticomplete to $\{x_1, \ldots ,x_q \} \backslash \{x_j\}$. Furthermore, for all distinct leaves $\ell, \ell'$ of $N$, $X_\ell$ is anticomplete to $X_{\ell'}$.

Let $I$ be the set of $i \in \{1, \ldots , k'\}$ such that $X_\ell \subseteq  C_i$ for some leaf $\ell$ of $N$. For each $i\in I$, let $\ell_i$ be the leaf of $N$ with $X_{\ell_i} \subseteq C_i$ (note that such a leaf is unique since $(X_v : v\in V(N))$ is $\{C_1,\ldots,C_{k'}\}$-spread).
For each $i\in I$, we define $x^i$ to be the vertex $x_j$ where $j \in \{1, \ldots ,q\}$ is chosen such that $\ell_i$ is adjacent to $n_j$ in $N$.
Thus there may be distinct values $i, i'\in I$ with $x^i =x^{i'}$.

Let $i \in I$. Since $X_{\ell_i}$ is $(\kappa + \epsilon)$-dominant and $(G,\mu)$ is $(\epsilon,5k')$-coherent, $N_G[X_{\ell_i}]\backslash N_G^{5k'}[x_1]$ is non-empty.
Therefore $X_{\ell_i} \backslash N_G^{5k' -1}[x_1]$ is non-empty.
Since $X_{\ell_i}$ contains a neighbour of $x^i$, the induced subgraph $G[X_{\ell_i} \cup \{x^i\}]$ is connected.
Observe that $q+4i \le 5k'$ and therefore $x^i \not\in N_G^{q+4i-1}[x_1]$ and $X_{\ell_i} \backslash N_G^{q+4i-1}[x_1]$ is non-empty.
Thus, there is a minimal path $P_i$ in $G[X_{\ell_i} \cup \{x^i\}]$ from $x^i$ to some vertex $u_i$ in $X_{\ell_i} \backslash N_G^{q+4i-1}[x_1]$.
By minimality, the distance between $u_i$ and $x_1$ is exactly $q+4i$.
Choose a vertex $b_i \in B_i$ adjacent to $u_i$.
Note that the distance between $x_1$ and $b_i$ in $G$ is between $q+4i-1$ and $q+4i+1$, so $\{b_i : i\in I\}$ is an independent set.
Let $c_i$ be the second vertex of $P_i$, that is, the vertex adjacent to $x^i$, and let $Q_i$ be a shortest path of $G[V(P_i) \cup \{b_i\} \backslash \{x^i\}]$ between $c_i$ and $b_i$. Thus $Q_i\backslash \{b_i\}$ is a path of $G[X_{\ell_i}]$. The subgraph of $G$ induced on $\{x_1, \ldots,x_q\}\cup  \{c_i :i\in I\}$ is a caterpillar isomorphic to $N$.

Next we will show that the subgraph $N'$ of $G$ induced on $\{x_1, \ldots,x_q\} \cup \bigcup_{i\in I} V(Q_i)$ is a isomorphic to a subdivision of $N$.
Let $i \in I$.
It is enough to show that $V(Q_i)\backslash \{c_i\}$ is anticomplete to $\{x_1, \ldots,x_q\}$, and for every $j\in I$ with $i<j$, $V(Q_i)$ is anticomplete to $V(Q_j)$.
Since $X_{\ell_i}$ is anticomplete to $\{x_1, \ldots ,x_q\} \backslash \{x^i\}$, and $V(Q_i) \backslash \{b_i\} \subseteq X_{\ell_i}$, we have that $V(Q_i) \backslash \{b_i\}$ is anticomplete to $\{x_1, \ldots ,x_q\} \backslash \{x^i\}$. By the minimality of $P_i$, $x^i$ has no neighbour in $V(Q_i) \backslash \{b_i\}$.
Furthermore, $b_i$ has no neighbour in $\{x_1, \ldots,x_q\}$, since the distance in $G$ between $x_1$ and $b_i$ is at least $q+4i-1>q$, and $G[\{x_1, \ldots,x_q\}]$ is a path. So, $V(Q_i)\backslash \{c_i\}$ is anticomplete to $\{x_1, \ldots,x_q\}$.
If $i$ is not the maximum element of $I$, then let $j\in I$ be such that $i<j$.
Certainly $V(Q_i)\backslash \{b_i\}$ and $V(Q_j)\backslash\{b_j\}$ are anticomplete, since $V(Q_i)\backslash\{b_i\}\subseteq X_{\ell_i}$, $V(Q_j)\backslash\{b_j\}\subseteq X_{\ell_j}$, and $X_{\ell_i}$ and $X_{\ell_j}$ are anticomplete. Also $b_i$ has no neighbour in $V(Q_j)\backslash\{b_j\}$ since the $k$-ladder is half-cleaned;
so it remains to check that $b_j$ has no neighbour in $V(Q_i)$.
By the minimality of $P_i$, the distance in $G$ between $x_1$ and each vertex of $V(Q_i)\backslash \{b_i\}$ is at most $q+4i$. So, the distance in $G$ between $x_1$ and each vertex of $V(Q_i)$ is at most $q+4i+1$. Since, the distance in $G$ between $x_1$ and $b_j$ is at least $q+4j-1> q+4i+2$, it follows that $b_j$ has no neighbour in $V(Q_i)$.
Thus, $V(Q_i)$ and $V(Q_j)$ are anticomplete, and it follows that the subgraph $N'$ of $G$ induced on $\{x_1, \ldots,x_q\} \cup \bigcup_{i\in I} V(Q_i)$ is a isomorphic to a subdivision of $N$.

For $1\le i \le k'$, let $C'_i =C_i\backslash N_G^{q+4k'+4}[x_1]$.
Since $k'\ge 2^{q+2}$, it follows that $5k'\ge q+4k'+4$, and so $(G,\mu)$ is $(\epsilon,q+4k'  + 4)$-coherent. Hence $\mu(N_G^{q+4k'+4}[x_1]) < \epsilon$. Thus for $1 \le i \le k'$, $\mu(C'_i) \ge \mu(C_i)-\epsilon \ge \kappa$.
Since every vertex in $V(N')$ has distance at most $q+4k'+1$ from $x_1$ in $G$, it follows that every vertex in $C'_i$ has distance at least three from $V(N')$ in $G$.
Let $B'_i$ be the set of vertices in $B_i \backslash \{b_i\}$ with no neighbours in $V(N')$. Since every vertex in $C'_i$ has a neighbour in $B_i$ and has distance at least three from $V(N')$ in $G$, it follows that $B_i'$ covers $C'_i$.

The sets
\[
A_i \ (i \in I), B_i' \ (i \in I), C_i' \ (i \in I)
\]
now form the desired $k$-ladder.
\end{proof}

For $I\subseteq [k]$, an \emph{$I$-tick} of a $k$-ladder $(A,B,C)$ is a collection of induced paths $(Q_i: i\in I)$ of $G$, each with a distinct end $q_i$, and a common end $x_I$ such that
\begin{itemize}
	\item the sets $(Q_i \backslash \{x_I\} : i \in I)$ are pairwise anticomplete,
    \item for each $i\in I$, $Q_i\backslash \{q_i\}$ is anticomplete to the $k$-ladder $(A,B,C)$,
    \item for each $i\in I$, $\{q_i\}$ is disjoint to the $k$-ladder $(A,B,C)$, and $N(q_i)\cap (A,B,C)\subseteq A_i$, and
    \item for each $i\in I$, $A_i\cup \{q_i\}$ is connected.
\end{itemize}

We call $x_I$ the \emph{centre} of the tick.
Given a tick $(Q_i: i\in I)$ of a $k$-ladder $(A,B,C)$ with centre $x_I$, for each $i\in I$ and $c\in C_i$, let $P_c$ be a path in $G[V(Q_i) \cup A_i \cup B_i \cup C_i]$ between $c$ and $x_I$ which is minimal subject to $V(P_c)\cap C_i = \{c\}$ and $|V(P_c)\cap B_i| = 1$. If for every $i\in I$ and $c\in C_i$, the path $P_c$ has odd length, then we say that $(Q_i: i\in I)$ is an \emph{odd tick} of the $k$-ladder $(A,B,C)$.
Similarly, if for every $i\in I$ and $c\in C_i$, the path $P_c$ has even length, then we say that $(Q_i: i\in I)$ is an \emph{even tick} of the $k$-ladder $(A,B,C)$.

Now, we derive the following from Lemma~\ref{EH:lem:caterpillarladder}.

\begin{lemma}\label{ticks}
Let $\epsilon,\kappa >0$, such that
$(2^{4pq+4q-1}-1)2^{4pq+4q-1}(2^{2^{4pq+4q-1}}(2\kappa + 2\epsilon) + \epsilon)\le 1$ and
$(2^{4pq+4q-1}-1)(2^{2^{4pq+4q-1}}(2\kappa + 2\epsilon) + \epsilon)+4\epsilon\le 1$.
Let $(G,\mu)$ be an $(\epsilon,5\cdot 2^{4pq+4q-1})$-coherent massed graph. Then there is a half-cleaned $pq$-ladder
\[
A_1, \ldots , A_{pq}, B_1, \ldots , B_{pq}, C_1, \ldots , C_{pq}
\]
in $G$, with pairwise anticomplete ticks $(Q_i^1: i\in I_1), \ldots, (Q_i^q: i\in I_q)$, such that
\begin{itemize}
    \item $\mu(C_i)\ge \kappa$ for $1\le i \le pq$,
    \item $I_1 \cup \cdots \cup I_q=\{1,\ldots , pq\}$,
    \item $|I_j|=p$ for each $1\le j \le q$, and
    \item either every tick is even, or every tick is odd.
\end{itemize}
\end{lemma}

\begin{proof}
    First observe that a forest consisting of $2q$ stars, each with $2p$ leaves is an induced subgraph of a caterpillar with $4pq+4q-1$ vertices, and $4pq$ leaves.
    So, by Lemma~\ref{EH:lem:caterpillarladder}, $(G,\mu)$ contains a half-cleaned $4pq$-ladder
    \[
    A_1,  \ldots , A_{4pq}, B_1, \ldots , B_{4pq}, C_1, \ldots , C_{4pq},
    \]
    with pairwise anticomplete ticks $(Q_i^1: i\in I_1), \ldots, (Q_i^{2q}: i\in I_{2q})$, such that
\begin{itemize}
    \item $\mu(C_i)\ge 2\kappa$ for $1\le i \le 4pq$,
    \item $I_1 \cup \cdots \cup I_{2q}=\{1,\ldots , 4pq\}$, and
    \item $|I_j|=2p$ for each $1\le j \le 2q$.
\end{itemize}
    
    Let $x_1,\ldots , x_{2q}$ be the centres of the ticks $(Q_i^1: i\in I_1), \ldots, (Q_i^{2q}: i\in I_{2q})$. For each $j\in [2q]$, $i\in I_j$, and $c\in C_{i}$, let $P_c$ be a path in $G[V(Q_i^j) \cup A_i \cup B_i \cup C_i]$ between $c$ and $x_j$ which is minimal subject to $V(P_c)\cap C_i = \{c\}$ and $|V(P_c)\cap B_i| = 1$.
    For every $i\in [4pq]$, since $\mu(C_i)\ge 2 \kappa$, there exists $C_i'\subseteq C_i$ such that $\mu(C_i')\ge \mu(C_i)/2 \ge \kappa$, and either for every $c\in C'_i$, the path $P_c$ has odd length, or for every $c\in C'_i$, the path $P_c$ has even length.
    Let $C'=\bigcup_{i\in [4pq]} C_i'$.
    
    Similarly, for every $j \in [2q]$, there exists a $I_j'\subseteq I_j$ with $|I_j'|=p$ such that either for every $i\in I_j'$, and $c\in C_{i}'$, the path $P_c$ has odd length, or for every $i\in I_j'$, and $c\in C_{i}'$, the path $P_c$ has even length. In other words, for every $j \in [2q]$, there exists a $I_j'\subseteq I_j$ with $|I_j'|=p$ such that the tick $(Q_i^j : i\in I_j')$ is either odd or even.

    Going one step further, there exists a $Q\subseteq [2q]$ with $|Q|=q$ such that either for every $j\in Q$, $(Q_i^{j}: i\in I_{j}')$ is an odd tick of $(A,B,C')$, or for every $j\in Q$, $(Q_i^{j}: i\in I_{j}')$ is an even tick of $(A,B,C')$.
    This gives us the desired ticks and $pq$-ladder.
\end{proof}

For a set $V$ and a non-negative integer $k$, we let ${V \choose k}$ denote the set of all $k$-element subsets of $V$.
We are now ready to prove Lemma~\ref{EH:lem:rcoherent}.

\begin{proof}[Proof of Lemma~\ref{EH:lem:rcoherent}.]
    Clearly it is enough to prove the lemma in the case that $H$ is a complete graph. So, let $H=K_n$.
    We may assume that $n\ge 2$.
    Let $\epsilon = 2^{-2^{4n^2}-8n^2}  n^{-2}$, $r= 5 \cdot 2^{4n^2}$, and let $(G,\mu)$ be an $(\epsilon , r)$-coherent massed graph.

    By Lemma~\ref{ticks}, there is a half-cleaned $n(n-1)$-ladder
\[
A_1, \ldots , A_{n(n-1)}, B_1, \ldots , B_{n(n-1)}, C_1, \ldots , C_{n(n-1)}
\]
in $G$, with pairwise anticomplete ticks $(Q_i^1: i\in I_1), \ldots, (Q_i^n: i\in I_n)$, such that
\begin{itemize}
    \item $\mu(C_i)\ge n(n-1)\epsilon$ for $1\le i \le n(n-1)$,
    \item $I_1 \cup \cdots \cup I_n=\{1,\ldots , n(n-1)\}$,
    \item $|I_j|=n-1$ for each $1\le j \le n$, and
    \item either every tick is even, or every tick is odd.
\end{itemize}

Let $x_1,\ldots , x_{n}$ be the centres of the ticks $(Q_i^1: i\in I_1), \ldots, (Q_i^{n}: i\in I_{n})$. For each $j\in [n]$, $i\in I_j$, and $c\in C_{i}$, let $P_c$ be a path in $G[V(Q_i^j) \cup A_i \cup B_i \cup C_i]$ between $c$ and $x_I$ which is minimal subject to $V(P_c)\cap C_i = \{c\}$ and $|V(P_c)\cap B_i| = 1$.
Note that either all such paths $P_c$ have even length, or they all have odd length. 
		For each such $j\in [n]$, $i\in I_j$, and $c\in C_{i}$, let $b_c$ be the unique vertex of $V(P_c)\cap B_i$.
		Observe that for each $j,j'\in [n]$ with $j\not=j'$, $i\in I_j$, $i'\in I_{j'}$ with $i<i'$, $c\in C_{i}$ and $c'\in C_{i'}$ with $c$ adjacent to $c'$, the subgraph $G[V(P_{c}) \cup V(P_{c'})] \backslash \{b_cb_{c'},cb_{c'}\}$ is an odd path between $x_j$ and $x_{j'}$, and so in particular, one of $G[V(P_{c}) \cup V(P_{c'})]$ or $G[(V(P_{c}) \cup V(P_{c'})) \backslash \{c, c'\}]$ is a fuzzy odd path between $x_j$ and $x_{j'}$.
		
		Arbitrarily order the elements ${[n] \choose 2}$ by $p_1,\ldots , p_{\frac{1}{2}n(n-1)}$.
		For $1 \le m \le \frac{1}{2}n(n-1)$, let $s_{m} \in I_i$, $s_{m}'\in I_j$ be such that $\bigcup_{\{i,j\}\in {[n] \choose 2}} \{s_{m}, s_{m}'\} = [n(n-1)]$.
		For $1 \le m \le \frac{1}{2}n(n-1)$, we will inductively find edges $c_{1}c_{1}',\ldots , c_{m}c_{m}'$ such that
		\begin{itemize}
			\item for every $1\le i \le m$, $c_{i} \in C_{s_i}$, and $c_{i}' \in C_{s_i'}$, and
			\item for every pair $1\le i < j \le m$, the distance in $G$ between $\{c_{i},c_{i}'\}$ and $\{c_{j},c_{j}'\}$ is at least 4.
		\end{itemize}
		
		Let $1\le m \le \frac{1}{2}n(n-1)$, and suppose we have already found such edges $c_{1}c_{1}',\ldots , c_{m-1}c_{{m-1}}'$.
		Since $(G, \mu)$ is $(\epsilon,4)$-coherent, for $s\in \{s_m,s_m'\}$ we have that
		\[
		\mu\left(C_{s} \backslash \bigcup_{i=1}^{m-1} N_G^4[\{c_{i}, c_{{i}}'\}]\right) \ge \mu(C_{s}) - \mu\left(\bigcup_{i=1}^{m-1} N_G^4[\{c_{{i}}, c_{{i}}'\}]\right)
		> n(n-1)\epsilon - 2(m-1)\epsilon > \epsilon.
		\]
		Therefore, since $(G,\mu)$ is $\epsilon$-coherent, there exists  vertices $c_{m} \in C_{s_m}\backslash \bigcup_{i=1}^{m-1} N_G^4[\{c_{{i}}, c_{{i}}'\}]$, and $c_{m}' \in C_{s_m'}\backslash \bigcup_{i=1}^{m-1} N_G^4[\{c_{{i}}, c_{{i}}'\}]$, with $c_{m}$ and $c_{m}'$ adjacent. Clearly the distance in $G$ between $\{c_{m}, c_{m}' \}$ and $\bigcup_{i=1}^{m-1} \{c_{{i}},c_{{i}}'\} $ is at least 5. Hence $c_{1}c_{1}',\ldots , c_{m}c_{m}'$ provides the desired collection of edges.
		So, by induction, we can find such edges $c_{1}c_{1}',\ldots , c_{{\frac{1}{2}n(n-1)}}c_{{\frac{1}{2}n(n-1)}}'$.
		
		For each $p_a=\{i,j\}\in {[n] \choose 2}$, observe that one of $G[V(P_{c_{a}}) \cup V(P_{c_a'})]$ or $G[(V(P_{c_{a}}) \cup V(P_{c_a'}))\backslash \{c_a,c_a'\}]$ is a fuzzy odd path between $x_i$ and $x_j$.
		So, for each $p_a=\{i,j\}\in {[n] \choose 2}$, choose $Q_{\{i,j\}}$ from $G[V(P_{c_{a}}) \cup V(P_{c_a'})]$ and $G[(V(P_{c_{a}}) \cup V(P_{c_a'}))\backslash \{c_a,c_a'\}]$ so that $Q_{\{i,j\}}$ is a fuzzy odd path between $x_i$ and $x_j$.
		Furthermore, for $\{i,j\}, \{i',j'\}\in {[n] \choose 2}$ with $p_a=\{i,j\}$, $p_{a'}=\{i',j'\}$, and $\{i,j\} \not= \{i',j'\}$, since the distance in $G$ between $\{c_{{a}}, c_{a}'\}$ and $\{c_{a'}, c_{a'}'\}$ is at least 5, we have that $V(Q_{\{i,j\}}) \backslash \{x_i,x_j\}$ is anticomplete to $V(Q_{\{i',j'\}}) \backslash \{x_{i'},x_{j'}\}$.
		Therefore, $G[\bigcup_{\{i,j\}\in {[n] \choose 2}} V(Q_{\{i,j\}}) ]$ is a proper fuzzy odd subdivision of $K_n$ as desired.
\end{proof}

\section{Focused massed graphs}\label{EH:sec:focus}

This section is dedicated to proving Lemma~\ref{EH:lem:focussedcoherent}.

Let $T'$ be a caterpillar and let $(G, \mu)$ be a massed graph containing an induced subgraph $T$ isomorphic to $T'$. Let $x_1,\ldots , x_q$ be the vertices of $T$ with degree more than one, and for each leaf $v$ of $T$, let $x^v$ be its neighbour in $\{x_1,\ldots, x_q\}$. Let $\kappa > 0$, and let $r\ge 1$ be an integer. If for each leaf $v$ of $T$ there is a subset $X_v \subseteq V(G)$ with the follow properties;
\begin{itemize}
    \item for each leaf $v\in L(T)$, $v\in X_v$, $x^v\not\in X_v$, and $X_v$ is anticomplete to $\{x_1, \ldots ,x_q\} \backslash \{x^v\}$,

    \item for all distinct leaves $u,v\in L(T)$, $X_u$ is anticomplete to $X_v$,

    \item for each leaf $v\in L(T)$, $x^v$ is an $r$-centre for $X_v \cup \{x^v\}$,

    \item for each leaf $v\in L(T)$, $v$ is the unique neighbour of $x^v$ in $X_v$, and

    \item for each leaf $v\in L(T)$, $X_v \cup \{x^v\}$ is $\kappa$-dominant,
\end{itemize}
then we say that $(T,(X_v)_{v\in L(T)})$ is a \emph{$(T', r, \kappa)$-frame} of $(G, \mu)$.

A result very similar to the following lemma is implicit in the proof of {\cite[7.2]{chudnovsky2021pure}} due to Chudnovsky, Scott, Seymour, and Spirkl. Since it is not given explicitly, we give a proof closely following {\cite[7.2]{chudnovsky2021pure}}.

\begin{lemma}\label{EH:lem:frame}
    Let $T'$ be a caterpillar with $t\ge 3$ vertices, and let $r \ge 1$ be an integer.
    Let $\epsilon, \delta, \kappa >0$ with $t\epsilon \le \delta \le \kappa \le 2^{-2^t}t^{-2^t}(2^{-2^t} (2^{-t} -\epsilon) -\epsilon )$.
    Then, every $(\delta,r)$-focused $\epsilon$-coherent massed graph $(G,\mu)$ contains a $(T', 3r+2, \kappa)$-frame.
\end{lemma}

\begin{proof}
    We may assume that $T'$ is rooted, and its head is an end of the path obtained from $T$ by deleting its leaves. Let $k= 2^t$, and let $\lambda = 1/k- \epsilon$.

    For each $1\le j \le k$ in order, choose some $Y_j\subset V(G) \backslash \left( \bigcup_{i<j} Y_i \right)$ minimal, subject to $\mu(Y_j) \ge \lambda$ (so with $\mu(Y_j)< \lambda + \epsilon)$.
    This is possible since
    \begin{align*}
    \mu\left(V(G) \backslash \left( \bigcup_{i<j} Y_i \right)\right)
    &\ge 1 - \mu\left( \bigcup_{i<j} Y_i \right) \\
    &\ge 1- \sum_{i<j} \mu(Y_i) \\
    &> 1- (j-1)(\lambda + \epsilon) \\
    &= 1-(j-1)/k \\
    &\ge 1/k \\
    &>\lambda.
    \end{align*}
    So, $Y_1,\ldots , Y_k$ are disjoint, and $\mu(Y_j) \ge \lambda$ for all $1\le j \le k$.

    For $0\le i \le k$, let $\kappa_i=2^{-i}t^{-i}(2^{-k}\lambda - \epsilon)$.
    Then, $\kappa_0 \ge \kappa_k = 2^{-2^t}t^{-2^t}(2^{-2^t} (2^{-t} -\epsilon) -\epsilon ) \ge \kappa \ge \delta$.
    As $\kappa_0 \ge \delta$ and $G$ is $(\delta,r)$-focused, $G$ is also $(\kappa_0,r)$-focused.
    So by Lemma~\ref{EH:lem:realization}, there is a $\{Y_1,\ldots ,Y_k\}$-spread $\kappa_0$-realization $(X_v : v\in V(T'))$ of $T'$ in $G$ such that $X_v$ has an $r$-centre for each $v\in V(T')$, except the head.
    Let $t_1,\ldots ,t_q$ be the vertices of degree more than one in $T'$, where $t_1$ is the head, and $t_{i}t_{i+1}$ is an edge for each $1\le i \le q-1$.
    Choose $x_1\in X_{t_1}$, and inductively for each $2\le i \le q$, choose $x_{i}\in X_{t_i}$ such that $x_i$ is adjacent to $x_{i-1}$.
    Then, $x_1,\ldots , x_q$ are the vertices of an induced path in $G$.

    Let $v_1,\ldots , v_\ell$ be the leaves of $T'$, and for each $v\in L(T')$, let $x^v$ be the vertex $x_j\in \{x_1,\ldots,x_q\}$ such that $v$ is adjacent to $t_j$ in $T'$.
    Since $X_v$ has an $r$-centre and $x^v$ has a neighbour in $X_v$, it follows that $x^v$ is a $(2r+1)$-centre of $X_v\cup \{x^v\}$.
    
    For each $v\in L(T')$, let $X^0_v=X_v\cup \{x^v\}$.
    For each $0\le i \le \ell$ in order, we shall inductively define $X^i_v$ for each $v\in L(T')$ satisfying the following:
    \begin{itemize}
        \item for each $v\in L(T')$, $x^v\in X_v^i$ and $X_v^i\backslash \{x^v\}$ is anticomplete to $\{x_1 , \ldots , x_q\} \backslash \{x^v\}$,

        \item for all distinct $u,v\in L(T')$, $X_u^i\backslash \{x^u\}$ is anticomplete to $X_v^i\backslash \{x^v\}$,

        \item for $1\le j \le i$, $x^{v_j}$ is a $(3r+2)$-centre of $X^i_{v_j}$, and $x^{v_j}$ has a unique neighbour in $X^i_{v_j}$,

        \item for $i < j \le \ell$, $x^{v_j}$ is a $(2r+1)$-centre of $X^i_{v_j}$, and

        \item for each $v\in L(T')$, $X^i_v$ is $\kappa_i$-dominant.
    \end{itemize}
	The sets $(X^0_v : v\in L(T'))$ satisfy these conditions.
    So, suppose that $1\le i \le \ell$, and we have defined $X^{i-1}_v$ for each $v \in V(T')$ as above.
    
    Let $v_j \in L(T') \backslash \{v_i\} = \{v_1,\ldots , v_\ell\} \backslash \{v_i\}$. If $j < i$, then let $r_j = 3r+ 2$, otherwise if $j > i$, then let $r_j = 2r+ 1$. So in either case $x^{v_j}$ is an $r_j$-centre of $X^{i-1}_{v_j}$. Choose $X^i_{v_j} \subseteq X^{i-1}_{v_j}$, minimal such that $X^i_{v_j}$ is $\kappa_i$-dominant and $x^{v_j}$ is an $r_j$-centre for $X^i_{v_j}$. It follows from the minimality that $\mu(N_G[X^i_{v_j}]) \le \kappa_i + 2\epsilon$ (since $\mu(N_G[v]) \le \mu(v) + \mu(N_G(v)) \le 2 \epsilon$ for every vertex $v$ of $G$).
    
    Since $X^{i-1}_{v_i}$ is $\kappa_{i-1}$-dominant, we have that $\mu(N_G[X^{i-1}_{v_i}])\ge \kappa_{i-1}$.
    Let 
    \[
    Y'=N_G[X^{i-1}_{v_i}] \backslash \left( N_G[\{x_1,\ldots , x_q\}] \cup \bigcup_{u\in L(T')\backslash \{v_i\}} N_G[X^i_u] \right).
    \]
    Then,
    \begin{align*}
    \mu(Y') &\ge \mu(N_G[X^{i-1}_{v_i}]) - \mu( N_G[\{x_1,\ldots , x_q\}] ) - \left( \sum_{ {u\in L(T')\backslash \{v_i\}} } \mu( N_G[X^i_u] ) \right) \\
    &\ge \kappa_{i-1} - 2 q\epsilon - (\ell - 1)(\kappa_i + 2 \epsilon) \\
    &=   2t\kappa_i  - 2 q\epsilon - (\ell - 1)(\kappa_i + 2 \epsilon) \\
    &=   (2t-\ell +1) \kappa_i  - 2 (q+\ell -1) \epsilon \\
    &\ge 4\kappa_i - 2 t\epsilon \\
    &\ge 2\kappa_i \\
    &\ge 2\kappa \\
    &\ge \delta.
    \end{align*}
    Since $G$ is $(\delta,r)$-focused, there is a subset $Y \subseteq Y'$ with $\mu(Y)\ge \mu(Y')/2 \ge \kappa_i$, such that $Y$ has an $r$-centre $y$ say.
    Since $x^{v_i}$ is a $(2r+1)$-centre for $X^{i-1}_{v_i}$, and $y\in N_G[X^{i-1}_{v_i}]$, there is an induced path $P$ of $G[X^{i-1}_v \cup \{y\}]$ between $x^{v_i}$ and $y$, of length at most $2r + 2$. Let $X^i_{v_i} = Y \cup V(P)$. Then $x^{v_i}$ is a $(3r + 2)$-centre for $X^i_{v_i}$.
    Since $x^{v_i}$ is anticomplete to $Y$ and has only one neighbour in $V(P)$, it follows that $x^{v_i}$ has only one neighbour in $X^i_{v_i}$. 
    Moreover, $X^i_{v_i}$ is $\kappa_i$-dominant since $\mu(N_G[X^i_{v_i}]) \ge \mu(X^i_{v_i}) \ge \mu(Y) \ge \kappa_i$. This completes the inductive definition.

    For each $v \in  L(T')$, let $c_v$ be the unique neighbour of $x^v$ in $X^\ell_v$.
    The subgraph $T$ induced on $\{x_1, \ldots ,x_q\} \cup \{c_v : v \in L(T')\}$ is isomorphic to $T'$.
    For each $c_v\in L(T)$, let $X_{c_v}= X_v^\ell$.
    Then, $(T,(X_{c_v})_{c_v\in L(T)})$ is a \emph{$(T', 3r+2, \kappa)$-frame} of $G$, as desired.
\end{proof}

Before proving Lemma~\ref{EH:lem:focussedcoherent}, we require one final lemma that allows us to find certain induced fuzzy odd paths.

\begin{lemma}\label{EH:lem:fuzz}
	Let $r\ge 1$ be an integer, and let $\epsilon>0$.
	Let $(G,\mu)$ be an $\epsilon$-coherent massed graph with $u\in X_u \subseteq A \subseteq V(G)$ and $v\in X_v \subseteq B \subseteq V(G)$ such that
	\begin{itemize}
		\item $u$ is an $r$-centre of $X_u$ and $v$ is an $r$-centre of $X_v$,
		
		\item $X_u$ is anti-complete to $X_v$,
		
		\item $A\subseteq N[X_u]$ and $B\subseteq N[X_v]$, and
		
		\item $\mu(A),\mu(B) \ge (r(r+1)(4^{r+1}+3) +1)\epsilon$.
	\end{itemize}
	Then $G$ contains an induced fuzzy odd path $P$ between $u$ and $v$ of length at most $2r+3$, such that $V(P)\subseteq A \cup B$.
\end{lemma}

\begin{proof}
	Since $(G,\mu)$ is $\epsilon$-coherent and $X_u$ is anti-complete to $X_v$, we may assume without loss of generality that $\mu(X_u) < \epsilon$.
	
	For $0\le j \le r$, let $L_A^j$ be the set of vertices in $X_u$ whose distance in $G[X_u]$ from $u$ is equal to $j$.
	For $0\le j \le r$, let $D_A^j = (A\backslash X_u) \cap (N(L_A^j) \backslash ( \bigcup_{0\le p < j} N(L_A^p))$.
	So, $\mu(D_A^0) + \cdots + \mu(D_A^{r}) \ge \mu(D_A^0 \cup \cdots \cup D_A^{r}) \ge \mu(A) - \mu(X_u) > (r(r+1)(4^{r+1}+3) +1)\epsilon - \epsilon = r(r+1)(4^{r+1}+3)\epsilon$.
	Choose $j_A$ minimum such that $\mu(D_A^{j_A})\ge (4^{r+1}+3)\epsilon$.
	Since $(G,\mu)$ is $\epsilon$-coherent, we have that $j_A\ge 1$.
	By the minimality of $j_A$, we have $\mu(\bigcup_{0\le p < j_A} D^p_A) \le \mu(D_A^0) + \cdots + \mu(D_A^{j_A -1}) < j_A (4^{r+1}+3)\epsilon \le r(4^{r+1}+3)\epsilon$.

	Let $B'= B \backslash (\bigcup_{0\le p < j_A} D^p_A)$, then
	\[
	\mu(B') \ge \mu(B) - \mu\left(\bigcup_{0\le p < j_A} D^p_A \right) 
	> (r(r+1)(4^{r+1}+3) +1)\epsilon  -  r(4^{r+1}+3)\epsilon
	\ge (3r+4)\epsilon.
	\]
	For $0\le j \le r$, let $L_B^j$ be the set of vertices in $X_v$ whose distance in $G[X_v]$ from $v$ is equal to $j$.
	For $0\le j \le r$, let $D_B^j = (B'\backslash X_v) \cap (N(L_B^j) \backslash ( \bigcup_{0\le p < j} N(L_B^p))$.
	
	Suppose that there exists a $0\le j_B \le r$ such that $\mu(N(L_B^{j_B})\cap D_A^{j_A}) \ge 4^{j_B}\epsilon$.
	Choose such a $j_B$ to be minimum, and let $Z= (N(L_B^{j_B})\cap D_A^{j_A}) \backslash \left( \bigcup_{0\le p < j_B} N(L_B^{p}) \right)$. Since $(G,\mu)$ is $\epsilon$-coherent, we have that $j_B\ge 1$.
	Then,
	\begin{align*}
	\mu(Z) &\ge \mu(N(L_B^{j_B})\cap D_A^{j_A}) - \left(  \sum_{p=0}^{j_B-1} \mu(N(L_B^{p})\cap D_A^{j_A})  \right) \\
	&> 4^{j_B}\epsilon - \sum_{p=0}^{j_B-1} 4^{p}\epsilon \\
	&= 4^{j_B}\epsilon - \frac{1}{3}(4^{j_B}-1)\epsilon \\
	&\ge 3\epsilon.
	\end{align*}
	Now, choose $L\subseteq L_B^{j_B}$ minimal such that $\mu(N(L) \cap Z) \ge \epsilon$.
	Then $\mu(N(L) \cap Z) < 2\epsilon$, and therefore $\mu(Z \backslash N(L)) \ge \mu(Z) - \mu(N(L) \cap Z) > 3\epsilon - 2\epsilon = \epsilon$.
	So, since $(G,\mu)$ is $\epsilon$-coherent, there exist adjacent vertices $b_1, b_2$, with $b_1\in Z \backslash N(L)$ and $b_2 \in N(L) \cap Z$.
	For $i\in \{1,2\}$, let $a_i$ be the vertex of $L_A^{j_A}$ adjacent to $b_i$, and let $P_i$ be a shortest path in $G[X_u]$ between $u$ and $a_i$.
	Then $P_1$ and $P_2$ both have length $j_A$.
	Let $c$ be a vertex of $L$ adjacent to $b_2$, and let $P_c$ be a shortest path in $G[X_v]$ between $v$ and $c$. 
	Then, $P_c$ has length $j_B$.
	Let $P_1'=G[V(P_1) \cup \{b_1,b_2\} \cup V(P_c)]$, and $P_2'=G[V(P_2) \cup \{b_2\} \cup V(P_c)]$. Clearly $V(P_1'),V(P_2') \subseteq A \cup B$. Also, $P_1' \backslash \{a_1b_2\}$ is a path between $u$ and $v$ of length $j_A+j_B+3\le 2r+3$ and $P_2'$ is a path between $u$ and $v$ of length $j_A+j_B+2\le 2r+2$.
	Therefore one of $P_1'$ or $P_2'$ provides the desired induced fuzzy odd path.
	
	So, we may now assume that $\mu(N(L_B^{j})\cap D_A^{j_A}) < 4^{j}\epsilon$ for each $0\le j \le r$.
	Therefore, $\mu(N(X_v)\cap D_A^{j_A} ) < \sum_{j=0}^{r} 4^{j}\epsilon = \frac{1}{3}(4^{r+1}-1)\epsilon <4^{r+1}\epsilon$. Let $Z= D_A^{j_A} \backslash N(X_v)$.
	Then $\mu(Z) \ge \mu(D_A^{j_A}) - \mu(N(X_v)\cap D_A^{j_A} ) \ge (4^{r+1}+3)\epsilon - 4^{r+1}\epsilon = 3\epsilon$.
	Since $Z$ is anticomplete to $X_v$, we have that $\mu(X_v) < \epsilon$.
	So, $\mu(D_B^0) + \cdots + \mu(D_B^{r}) \ge \mu(D_B^0 \cup \cdots \cup D_B^{r}) \ge \mu(B') - \mu(X_v) \ge (3r+4)\epsilon - \epsilon = 3(r+1)\epsilon$.
	Now, choose $0\le j_B \le r$ minimum such that $\mu(D_B^{j_B}) \ge 3\epsilon$.
    Since $\mu(Z) \ge 3\epsilon$, there exists a minimal subset $D$ of $D_B^{j_B}$ such that $\mu(N(D)\cap Z) \ge \epsilon$.
	Then $\mu(N(D)\cap Z)< 2\epsilon$, and therefore $\mu(Z \backslash N(D)) \ge \mu(Z) - \mu(N(D)\cap Z) > 3\epsilon -2\epsilon = \epsilon$.
	So, since $(G,\mu)$ is $\epsilon$-coherent, there exist adjacent vertices $b_1, b_2$, with $b_1\in Z \backslash N(D)$ and $b_2 \in N(D) \cap Z$.
	For $i\in \{1,2\}$, let $a_i$ be the vertex of $L_A^{j_A}$ adjacent to $b_i$, and let $P_i$ be a shortest path in $G[X_u]$ between $u$ and $a_i$.
	Then $P_1$ and $P_2$ both have length $j_A$.
	Let $c_1$ be a vertex of $D$ adjacent to $b_2$, and let $c_2$ be a vertex of $L_B^{j_B}$ adjacent to $c_2$. Let $P_c$ be a shortest path in $G[X_v]$ between $v$ and $c_2$. 
	Then, $P_c$ has length $j_B$.
	Let $P_1'=G[V(P_1) \cup \{b_1,b_2,c_1\} \cup V(P_c)]$, and $P_2'=G[V(P_2) \cup \{b_2,c_1\} \cup V(P_c)]$. Clearly $V(P_1'),V(P_2') \subseteq A \cup B$. Also, $P_1' \backslash \{a_1b_2\}$ is a path between $u$ and $v$ of length $j_A+j_B+4\le 2r+4$ and $P_2'$ is a path between $u$ and $v$ of length $j_A+j_B+3\le 2r+3$.
	Therefore one of $P_1'$ or $P_2'$ provides the desired induced fuzzy odd path.
\end{proof}

A \emph{star} is a tree of radius $1$. We call the non-leaf vertex of a star its \emph{centre}.
We are now ready to prove Lemma~\ref{EH:lem:focussedcoherent}. The rest of the proof adapts that of another lemma of Chudnovsky, Scott, Seymour, and Spirkl {\cite[7.1]{chudnovsky2021pure}}.

\begin{proof}[Proof of Lemma~\ref{EH:lem:focussedcoherent}.]
    Clearly it is enough to prove the lemma in the case that $H$ is a complete graph on $n\ge 3$ vertices. So, let $H=K_n$, and let $r\ge 1$ be an integer.
    Now let $\epsilon, \delta, \kappa >0$ be such that
    \begin{itemize}
        \item $(n^2+n-1)\epsilon \le \delta \le \kappa \le 2^{-2^{n^2+n-1}}(n^2+n-1)^{-2^{n^2+n-1}}(2^{-2^{n^2+n-1}} (2^{-{(n^2+n-1)}} -\epsilon) -\epsilon )$
        , and
    	\item $2^{\frac{1}{2}n(n-1)}n^{n(n-1)} \left( 6 n^2(r+3) \epsilon   
    	+ ((3r+2)(3r+3)(4^{3r+3}+3) +1)\epsilon \right) \le \kappa$.
    \end{itemize}
	Such $\epsilon, \delta, \kappa >0$ can clearly be chosen when $\epsilon$ is sufficiently small.
    Then, let $(G,\mu)$ be a $(\delta , r)$-focused $\epsilon$-coherent massed graph.
    
	For every $i\in [n]$, let $S_i'$ be a star with centre $w_i'$ and leaves $L(S_i')=\{u_{\{i,j\}}' : j\in [n] \backslash [i] \} \cup \{v_{\{j,i\}}' : j\in [i-1] \}$. Let $T'$ be a caterpillar with $n^2+n-1$ vertices that contains the disjoint union of $S_1',\ldots , S_n'$ as an induced subgraph (each of these $n$ stars has $n$ vertices, and then we need an extra $n-1$ vertices to join them together into a caterpillar).
	Then by Lemma~\ref{EH:lem:frame}, $(G, \mu)$ contains a $(T', 3r+2, \kappa)$-frame $(T,(X_v)_{v\in L(T)})$.
	For each $i\in [n]$, let $S_i$ be the induced star of $T$ corresponding to the induced star $S_i'$ of $T'$. Let $w_i$ be the centre of $S_i$ and let $L(S_i) = \{u_{\{i,j\}} : j\in [n] \backslash [i] \} \cup \{v_{\{j,i\}} : j\in [i-1] \}$ for each $i\in [n]$.
	For each $i\in [n]$, and each leaf $w$ of $S_i$, let $x^w=w_i$.

	Let $\kappa_0 =\kappa$.
	For each $1 \le i \le \frac{1}{2}n(n-1)$ in order, let
	\[
	\kappa_i = \frac{1}{n^2} \left( \kappa_{i-1} - 6 n^2(r+3) \epsilon   
	- ((3r+2)(3r+3)(4^{3r+3}+3) +1)\epsilon \right).
	\]
	For each $1 \le i \le \frac{1}{2}n(n-1)$, we have that $\kappa_{i} \ge \frac{1}{2n^2}\kappa_{i-1}$. So, $\kappa_0 > \cdots > \kappa_{\frac{1}{2}n(n-1)}
	\ge 6 n^2(r+3) \epsilon   
	+ ((3r+2)(3r+3)(4^{3r+3}+3) +1)\epsilon$.
	Arbitrarily order the elements of ${[n] \choose 2}$ by $p_1,\ldots , p_{\frac{1}{2}n(n-1)}$.
	For each $w \in L(T)$, let $X^0_w = X_w \cup \{x^w\}$.
	For $1 \le i \le \frac{1}{2}n(n-1)$ in order, we define $X^i_w$ $(w\in L(T) \backslash \{u_{p_1}, \ldots, u_{p_{i-1}}, v_{p_1}, \ldots , v_{p_{i-1}}  \})$ and $P_i$ inductively as follows.
	We assume $P_1, \ldots ,P_{i-1}$ and $X^{i-1}_w$ $(w\in L(T) \backslash \{u_{p_1}, \ldots, u_{p_{i-1}}, v_{p_1}, \ldots , v_{p_{i-1}}  \})$ have been defined, such that
	\begin{itemize}
		\item for $1 \le h \le i-1$, $P_h$ is an induced fuzzy odd path between $u_{p_h}$ and $v_{p_h}$, of length at most $6r+ 7$,
		
		\item for $1\le h \le i-1$, $V(P_h)\backslash \{u_{p_h},v_{p_h}\}$ is anticomplete to $V(T)\backslash \{u_{p_h},v_{p_h}\}$,
		
		\item for $1\le h < h' \le i-1$, $V(P_h)$ is anticomplete to $V(P_{h'})$,
		
		\item for $1\le h \le i-1$ and $v\in  \{u_{p_i}, \ldots, u_{p_{\frac{1}{2}n(n-1)}}, v_{p_i}, \ldots , v_{p_{\frac{1}{2}n(n-1)}}  \}$, $V(P_h)$ is anticomplete to $X^{i-1}_v$,
		
		\item for $w \in \{u_{p_i}, \ldots, u_{p_{\frac{1}{2}n(n-1)}}, v_{p_i}, \ldots , v_{p_{\frac{1}{2}n(n-1)}}  \}$, $X^{i-1}_w$ is $\kappa_{i-1}$-dominant, and
		
		\item for $w \in \{u_{p_i}, \ldots, u_{p_{\frac{1}{2}n(n-1)}}, v_{p_i}, \ldots , v_{p_{\frac{1}{2}n(n-1)}}  \}$, $\{x^w, w\}\subseteq X^{i-1}_w$ and $x^w$ is a $(3r+2)$-centre for $X^{i-1}_w$.
	\end{itemize}

    For each $w \in \{u_{p_{i+1}}, \ldots, u_{p_{\frac{1}{2}n(n-1)}}, v_{p_{i+1}}, \ldots , v_{p_{\frac{1}{2}n(n-1)}}  \}$, choose $X^i_w \subseteq X^{i-1}_w$, minimal such that $X^i_w$ is $\kappa_i$-dominant and $x^w$ is a $(3r+2)$-centre for $X^i_w$. By deleting a vertex in $X^i_w$ with maximum $G[X^i_w]$-distance from $x^w$, the minimality of $X^i_w$ implies that $\mu(N_G[X^i_w]) \le \kappa_i+ 2 \epsilon$.
    We also have that
    \begin{align*}
    \mu\left( N_G\left[ \{w_1, \ldots ,w_n\} \cup \bigcup_{h<i} V(P_h)  \right]  \right) 
    &\le 2 (n+ (i-1)(6r+8))\epsilon \\
    &\le  2 (n+ \frac{1}{2}n(n-1)(6r+8))\epsilon \\
    &<   6 n^2(r+2) \epsilon.
    \end{align*}
    Let $u=u_{p_i}$ and $v =v_{p_i}$. Then, let $Y$ be the set of all vertices at distance at most one from $\{w_1, \ldots ,w_n\} \cup \bigcup_{h<i} V(P_h)$ and $X^i_w$ for each $w \in \{u_{p_{i+1}}, \ldots, u_{p_{\frac{1}{2}n(n-1)}}, v_{p_{i+1}}, \ldots , v_{p_{\frac{1}{2}n(n-1)}}  \}$.
    Observe that,
    \[
    \mu(Y)\le 6 n^2(r+2) \epsilon  + n(n-1)(\kappa_i+ 2 \epsilon) <  6 n^2(r+3) \epsilon  + n^2\kappa_i.
    \]
    
    Let $A= N_G[X^{i-1}_u] \backslash Y$, and $B= N_G[X^{i-1}_v] \backslash Y$.
    Then, since $X^{i-1}_u$ is $\kappa_{i-1}$-dominant, it follows that
    \[
    \mu(A) \ge \mu(N_G[X^{i-1}_u]  ) - \mu(Y) > \kappa_{i-1} - 6 n^2(r+3) \epsilon  - n^2\kappa_i \ge ((3r+2)(3r+3)(4^{3r+3}+3) +1)\epsilon.
    \]
    Similarly, $\mu(B) > ((3r+2)(3r+3)(4^{3r+3}+3) +1)\epsilon$.
    So, by Lemma~\ref{EH:lem:fuzz} (with $r$ replaced by $3r+2$), there exists an induced fuzzy odd path $P_i$ between $u$ and $v$ of length at most $6r+7$ with $V(P_i) \subseteq A\cup B$.
    Since $V(P_i) \subseteq A\cup B \subseteq N_G[X_u^{i-1}] \cup N_G[X_v^{i-1}]$, $P_i$ is anticomplete to $V(T)\backslash \{u,v\}$, $V(P_h)$ for each $1\le h \le i-1$, and $X_w^i \subseteq X_w^{i-1}$ for each $w \in \{u_{p_i}, \ldots, u_{p_{\frac{1}{2}n(n-1)}}, v_{p_i}, \ldots , v_{p_{\frac{1}{2}n(n-1)}}  \}$. Therefore, $P_i$ provides the desired induced fuzzy odd path.
    This completes the inductive definition of $P_1, \ldots , P_{\frac{1}{2}n(n-1)}$.
    
    It now follows that $G[\{w_1,\ldots w_n\} \cup \bigcup_{j=1}^{\frac{1}{2}n(n-1)} V(P_j) ]$ is a proper fuzzy odd subdivision of $H$, as desired.
\end{proof}

This completes the proof of Theorem~\ref{mainfuzzymass}, and thus also of Corollary~\ref{EH:masscol}, Theorem~\ref{main}, Corollary~\ref{EH:main}, and Theorem~\ref{mainfuzzy}.

\section*{Acknowledgements}

The author thanks Jim Geelen and Rose McCarty for a number of helpful discussions on pivot-minors.
The author also thanks the anonymous referees for helpful comments, in particular for spotting a serious error in an earlier version of this manuscript.

\bibliographystyle{amsplain}

\end{document}